\documentclass[11pt]{article}
\usepackage{cite}
\usepackage{amsmath,amssymb,amsfonts}
\usepackage{algorithmic}
\usepackage{graphicx,lipsum}
\usepackage{url}
\usepackage{textcomp}
\usepackage[left=1in,top=.7in,right=1in,bottom=1in]{geometry}
\date{}

\DeclareMathOperator{\vect}{vec}
\DeclareMathOperator{\sgn}{sgn}

\def\bc{\mbox{\boldmath $c$}}
\def\bn{\mbox{\boldmath $n$}}
\def\bs{\mbox{\boldmath $s$}}
\def\bq{\mbox{\boldmath $q$}}
\def\bp{\mbox{\boldmath $p$}}
\def\ba{\mbox{\boldmath $a$}}
\def\bb{\mbox{\boldmath $b$}}
\def\bu{\mbox{\boldmath $u$}}

\def\br{\mbox{\boldmath $r$}}

\def\bv{\mbox{\boldmath $v$}}

\def\bw{\mbox{\boldmath $w$}}

\def\bt{\mbox{\boldmath $t$}}

\def\bc{\mbox{\boldmath $c$}}
\def\bd{\mbox{\boldmath $d$}}

\def\bB{\mbox{\boldmath $B$}}

\def\bQ{\mbox{\boldmath $Q$}}

\def\bK{\mbox{\boldmath $K$}}

\def\bA{\mbox{\boldmath $A$}}

\def\bvarphi{\mbox{\boldmath $\varphi$}}

\def\bbeta{\mbox{\boldmath $\beta$}}

\DeclareMathOperator{\sign}{sgn}

\def\bm{\mbox{\boldmath $m$}}

\def\bsigma{\mbox{\boldmath $\sigma$}}

\def\bxi{\mbox{\boldmath $\xi$}}

\def\bt{\mbox{\boldmath $t$}}
\def\bn{\mbox{\boldmath $n$}}

\def\bxi{\mbox{\boldmath $\xi$}}

\def\brho{\mbox{\boldmath $\rho$}}

\def\bC{\mbox{\boldmath $C$}}
\def\bA{\mbox{\boldmath $A$}}
\def\bB{\mbox{\boldmath $B$}}

\def\bI{\mbox{\boldmath $I$}}
\def\bt{\mbox{\boldmath $t$}}
\def\bn{\mbox{\boldmath $n$}}

\def\bxi{\mbox{\boldmath $\xi$}}

\def\bI{\mbox{\boldmath $I$}}

\def\bR{\mbox{\boldmath $R$}}

\def\bg{\mbox{\boldmath $g$}}
\def\ba{\mbox{\boldmath $a$}}

\def\bp{\mbox{\boldmath $p$}}

\def\bv{\mbox{\boldmath $v$}}

\def\bLambda{\mbox{\boldmath $\Lambda$}}

\newtheorem{proposition}{Proposition}

\def\BibTeX{{\rm B\kern-.05em{\sc i\kern-.025em b}\kern-.08em
    T\kern-.1667em\lower.7ex\hbox{E}\kern-.125emX}}

\begin{document}
\title{The analytical subtraction approach for solving the forward problem in EEG}
\author{L. Beltrachini
\thanks{L. Beltrachini is with the Cardiff University Brain Research Imaging Centre (CUBRIC), School of Physics and Astronomy, Cardiff University, Cardiff CF24 3AA, UK. (email: BeltrachiniL@cardiff.ac.uk)}}

\maketitle

\begin{abstract}


Objective: The subtraction approach is known for being a theoretically-rigorous and accurate technique for solving the forward problem in electroencephalography by means of the finite element method. One key aspect of this approach consists of computing integrals of singular kernels over the discretised domain, usually referred to as potential integrals. Several techniques have been proposed for dealing with such integrals, all of them approximating the results at the expense of reducing the accuracy of the solution. In this paper, we derive analytic formulas for the potential integrals, reducing approximation errors to a minimum.

Approach: Based on volume coordinates and Gauss theorems, we obtained parametric expressions for all the element matrices needed in the formulation assuming first order basis functions defined on a tetrahedral mesh. This included solving potential integrals over triangles and tetrahedra, for which we found compact and efficient formulas.

Main results: Comparison with numerical quadrature schemes allowed to test the advantages of the methodology proposed, which were found of great relevance for highly-eccentric sources, as those found in the somatosensory and visual cortices. Moreover, the availability of compact formulas allowed an efficient implementation of the technique, which resulted in similar computational cost than the simplest numerical scheme.

Significance: The analytical subtraction approach is the optimal subtraction-based methodology with regard to accuracy. The computational cost is similar to that obtained with the lowest order numerical integration scheme, making it a competitive option in the field. The technique is highly relevant for improving electromagnetic source imaging results utilising individualised head models and anisotropic electric conductivity fields without imposing impractical mesh requirements.

\end{abstract}


\section{Introduction}

The forward problem in electroencephalography (EEG-FP) is a cornerstone of quantitative brain activity characterisation based on electromagnetic recordings. It consists in computing the electric potential function in the head generated by a set of known current generators representing time-varying charge distributions in the cortex at the mesoscopic level. Mathematically, the EEG-FP is described by a Poisson-like equation (subject to a Neumann boundary condition) depending on a particular combination of source and conductor model parameters~\cite{Beltrachini2019b}. It is widely accepted in the literature to consider dipoles as the generators of EEG signals~\cite{deMunck1988a,Hamalainen1993,Baillet2001}. Under such assumption, analytical solutions to the EEG-FP can be obtained for multilayered spherical head models with piecewise anisotropic electrical conductivity~\cite{deMunck1993}. However, it is acknowledged that the adoption of individualised head representations confer noticeable advantages over spherical models in the analysis EEG (either invasive~\cite{vonEllenrieder2012} or non-invasive~\cite{Cuffin1990,Huiskamp1999,Vatta2010,Beltrachini2011}), leading to the use of numerical techniques. 

Several numerical methodologies were presented for solving the EEG-FP, from which the finite element method (FEM) stands out of the rest. The reason is based on its flexibility for incorporating arbitrary geometries and heterogeneous and anisotropic electrical conductivity fields, which were shown to impact the current flow through biological tissues~\cite{Wolters2006,vonEllenrieder2012,vonEllenrieder2014b}. One key aspect in the solution of the EEG-FP by means of the FEM (as well as with any other numerical method) is related to the singularity introduced by the dipolar source model. Within few existing options for handling such singularities, the subtraction approach is the only that allows the use of standard FE basis functions for simulating dipoles in arbitrary locations while guaranteeing the existence and uniqueness of the solution~\cite{Beltrachini2019a,Beltrachini2019b,Wolters2007}. Subtraction-based techniques were shown to provide highly accurate solutions to the EEG-FP, improving those obtained by the partial integration method, and comparable to using the Saint Venant's principle~\cite{Wolters2007c,Lew2009,Vorwerk2012}. Solving the EEG-FP by means of the subtraction FEM requires the solution of surface and volume integrals with singular integrands, usually referred to in electromagnetism as {\it potential integrals}~\cite{Wilton1984}. So far, researchers have managed to work with approximations of these integrals obtained either projecting the singular function~\cite{Wolters2007} or its gradient~\cite{Beltrachini2019a} on the FE space, or utilising numerical quadrature schemes~\cite{Drechsler2009,Beltrachini2019b}. However, such approximations are known for introducing unwanted errors into the solution, requiring high resolution meshes and/or prohibitive computational resources to achieve accurate results.

In this paper, we derive analytical expressions for all the integrals needed in the subtraction FE formulation of the EEG-FP utilising linear basis functions defined on a tetrahedral mesh. This is done by exploiting the volume coordinate system and both surface and volume Gauss theorems to reduce integrals to 1D, where they become analytically solvable. Simple and compact expressions are obtained for all the integrals involved in the computation of the stiffness matrix and source vectors, avoiding any error other than that given by the domain discretisation itself. This makes the proposed method, coined {\it analytical subtraction}~(AS) approach, the optimal with regard to accuracy. We illustrate the advantages of the AS technique over the use of numerical quadrature formulas at the element level, and demonstrate how these differences impact in the overall solution of the EEG-FP. Additionally, we show that the expressions allow to reduce the computational requirements compared to the use of a numerical quadrature scheme, making it a competitive option. 

The rest of the paper is organised as follows: in Section~\ref{sec:2} we review the subtraction formulation of the EEG-FP, and present the integrals needed in the FE discretisation. The solution of these integrals require different tools, and are therefore treated separately. Analytical expressions for the element stiffness matrix are derived in Section~\ref{sec:3}, whereas the corresponding to the source vector are obtained in Section~\ref{sec:4}. Local and global experiments are described in Section~\ref{sec:6}, and the corresponding results presented in Section~\ref{sec:7}. Finally, we discuss the relevance of the methodology in Section~\ref{sec:8}, and summarise the conclusions in Section~\ref{sec:9}.

\section{EEG forward problem} \label{sec:2}

\subsection{Differential formulation}
The EEG-FP consists of finding the electric potential function~$u(\br)$ due to a current source with density~$s(\br)$ defined over the domain~$\Omega$ (i.e. the head), with boundary~$\Gamma$. Let~$\overline{\bsigma}(\br)$ be the rank-2 conductivity tensor field within~$\Omega$, and~$\hat{\bn}(\br)$ the unitary vector normal to~$\Gamma$. Under generally accepted assumptions (as the {\it quasistatic} and the {\it point electrode model} approximations), the EEG-FP reduces to find~$u(\br)$ satisfying $\nabla \cdot \big(  \overline{\bsigma} (\br) \nabla u(\br) \big)=-s(\br)$ ($\br\in\Omega$) together with the boundary condition $\left\langle \overline{\bsigma} (\br) \nabla u(\br),\hat{\bn}(\br)  \right\rangle=0$ ($\br \in \Gamma$)~\cite{deMunck1993,Beltrachini2019b}. In the case of assuming a dipolar source located in~$\br_0$ with dipolar moment~$\bq$, $s(\br)=-\left\langle\bq,\nabla \delta(\br - \br_0)\right\rangle$.

The subtraction method consists of avoiding the singularity in $s(\br)$ by separating~$\Omega$ into two subsets, one surrounding the source, namely~$\Omega^\infty$, with homogeneous conductivity~$\overline{\bsigma}^\infty=\overline{\bsigma}(\br_0)$, and the other being $\Omega^c=\Omega\backslash \Omega^\infty$ (the complement of $\Omega^\infty$ in $\Omega$) with electrical conductivity~$\overline{\bsigma}^c(\br)=\overline{\bsigma}(\br)-\overline{\bsigma}^\infty$. This  allows to express the electric potential as the sum of two terms, $u(\br)=u^c(\br)+u^\infty(\br)$, where~$u^\infty(\br)$ is the singularity potential generated by a source in an unbounded homogeneous conductor with conductivity~$\overline{\bsigma}^\infty$ (for which analytical expressions exist), and~$u^c(\br)$ is the correction potential satisfying $\nabla\cdot \big(\overline{\bsigma}(\br) \nabla u^{c}(\br)\big)=-\nabla\cdot \big(\overline{\bsigma}^{c}(\br) \nabla u^{\infty}(\br)\big)$ ($\br\in\Omega$) and subject to $\left\langle\overline{\bsigma}(\br)\nabla u^{c}(\br),\hat{\bn}(\br)\right\rangle=-\left\langle\overline{\bsigma}(\br)\nabla u^{\infty}(\br),\hat{\bn}(\br)\right\rangle$ ($\br\in\Gamma$)~\cite{Wolters2007,Drechsler2009,Beltrachini2019b}. The problem then results in finding the correction potential by means of a numerical method, after which the singularity potential is added.

\subsection{FE discretisation}

The FE formulation of the EEG-FP relies on the variational form of the subtraction version. This can be obtained by multiplying the corresponding differential equation by a test function~$v$ belonging to a suitable space~$H$, and then integrating over~$\Omega$. After utilising the divergence theorem and the boundary condition~\cite{Wolters2007,Drechsler2009}, the variational formulation results in finding~$u^c(\br)\in H$ such that, for all~$v(\br)\in H$, satisfies $a\left(u^c,v \right)=l(v)$, where $a:H\times H\rightarrow \mathbb{R}$ is the bilinear form given by
\begin{equation}\label{eq_sesquilinear}
a\left( u,v\right)=\int_\Omega \left\langle  \overline{\bsigma}(\br) \nabla u(\br),\nabla v(\br) \right\rangle d\br,
\end{equation}
and $l:H\rightarrow \mathbb{R}$ is the linear form defined as
\begin{align}
l(v)=-\int_\Omega \left\langle  \overline{\bsigma}^c(\br) \nabla u^{\infty}(\br),\nabla v(\br) \right\rangle d\br \nonumber\\
-\int_\Gamma v(\br) \left\langle  \overline{\bsigma}^{\infty} \nabla u^{\infty}(\br),\hat{\bn}(\br)\right\rangle d\br. \label{eq_l}
\end{align} 

To proceed with the FEM, a discretisation~$\mathcal{T}$ of~$\Omega$ is needed. This discretisation is composed by a set of nodes~$\bp_i$ ($i=1,\dots,N$) and elements~$\mathcal{T}_j$ ($j=1,\dots,N_e$) defined upon these nodes. This tessellation allows to construct a discretised FE space~$V_N\subset H$ where to compute the numerical solution. More explicitly, we choose $V_N=span\{\varphi_i(\br): i=1,\dots,N \}$, with~$\varphi_i(\br)$ being piecewise functions such that~$\varphi_i(\bp_j)=\delta_{ij}$~\cite{Wolters2007}. Then, we look for~$\widetilde{u}^c(\br)\in V_N$ (an approximation of $u^c(\br)\in H$) satisfying, for all $v(\br)\in V_N$, $a(\widetilde{u}^c,v)=l(v)$. This leads to solve the system of linear equations
\begin{equation}\label{eq:system}
\bK \bu^c = \bb,
\end{equation}
where $\bK\in\mathbb{R}^{N\times N}$ is the {\it stiffness matrix} defined by $K_{ij}=a(\varphi_i(\br),\varphi_j(\br))$, $\bb \in \mathbb{R}^N$ is the {\it source vector} with elements $b_i=l(\varphi_i(\br))$, and~$\bu^c\in\mathbb{R}^N$ is the vector containing the numerical approximation of the correction potential on the mesh nodes (i.e. $u^c(\br)\approx \widetilde{u}^c (\br)=\sum_{i=1}^N \varphi_i(\br)u^c_i$). The most common scenario available in the literature (and the one adopted in this work) is to utilise linear basis functions defined on a tetrahedral discretisation. In this case, the arrays $\bK$ and~$\bb$ needed in~\eqref{eq:system} are obtained by assembling the corresponding element matrices, i.e. those found by integrating over individual tetrahedrons $\mathcal{T}_j$ ($j=1,\dots,N_e$) or surface triangles $T_k$ ($k=1,\dots,N_s$). For simplicity, we split the source vector into two terms representing the surface and volume integrals, i.e. $\bb=-(\bb_s+\bb_v)$, allowing their individual analysis.

Without loss of generality, we assume that $\Omega^\infty$ is an isotropic medium, i.e. $\overline{\bsigma}^\infty=\sigma^\infty \bI_3$, with~$\bI_n$ being the $n\times n$ identity matrix. Moreover, we consider dipolar source models to represent current generators, as is standard in the literature~\cite{Hamalainen1993,deMunck1988a}. In this case, $u^\infty(\br)=(4\pi \sigma^\infty)^{-1}f(\br)$, with $f(\br)= \bq \cdot \bR R^{-3}$, $\bR=\br-\br_0$, and $R=|\bR|$. Under these assumptions, the element matrices are
\begin{gather}
\bK_k^e=\int_{\mathcal{T}_k} \left\langle  \overline{\bsigma}(\br) \nabla \bvarphi(\br),\nabla \bvarphi(\br) \right\rangle d\br, \label{eq_ke}\\
\bb_k^v=\frac{1}{4\pi\sigma^\infty}\int_{\mathcal{T}_k} \left\langle  \overline{\bsigma}^c(\br) \nabla f(\br),\nabla \bvarphi(\br) \right\rangle d\br,\label{eq_bv}\\
\bb_j^s=\frac{1}{4\pi}\int_{T_j} \bvarphi(\br) \left\langle \nabla f(\br),\hat{\bn}(\br)\right\rangle d\br,\label{eq_bs}
\end{gather} 
where $k=1,\dots,N_e$, $j=1,\dots,N_s$, and $\bvarphi(\br)$ is the linear interpolation vector function with elements $\varphi_i(\br)$, $i=1,\dots,d_i$ [$d_i=4$ for~\eqref{eq_ke} and~\eqref{eq_bv}; $d_i=3$ for~\eqref{eq_bs}; see~\eqref{eq:phi}]. In the following sections, we derive analytical expressions for~\eqref{eq_ke}--\eqref{eq_bs} based on the application of linear transforms and Gauss theorems.

\section{Stiffness matrix}\label{sec:3}

The computation of the element stiffness matrix is easily performed by employing the {\it volume coordinate system}~\cite{Hutton2003,Beltrachini2015a}. Let $V$ be the volume of an arbitrary tetrahedron $\mathcal{T}$ with nodes $\bp_i$ ($i=1,\dots,4$). For any point $\br\in \mathcal{T}$, we define the volume coordinates $\xi_i=V_i/V$ ($i=1,\dots,4$), where $V_i$ is the volume of the tetrahedron with nodes $\br$ and $\bp_j$, with $j\neq i$ ($j=1,\dots,4$). These coordinates are given by 
\begin{equation*}
\xi_i=\frac{a_i+b_i x_1 +c_i x_2 +d_i x_3}{6V},\ i=1,2,3,4,
\end{equation*}
where $\br=[x_1,x_2,x_3]^T$ and $a_i$, $b_i$, $c_i$, and $d_i$ are real coefficients depending on $\bp_i$~\cite{Hutton2003}. In matrix form,
\begin{equation}\label{transf_mat}
\bxi=\frac{1}{6V}\big( \tilde{\ba}+\bLambda^T \br \big),
\end{equation}
where $\tilde{\ba}$, $\tilde{\bb}$, $\tilde{\bc}$, and $\tilde{\bd}$ are the $4\times 1$ vectors with elements $a_i$, $b_i$, $c_i$, and $d_i$, respectively ($i=1,\dots,4$), $\bLambda=[\tilde{\bb},\tilde{\bc},\tilde{\bd}]^T$, and $\bxi=[\xi_1,\xi_2,\xi_3,\xi_4]^T$. It is easy to note that~\eqref{transf_mat} transforms $\mathcal{T}$ into a normalised tetrahedron $\mathcal{T}_{n}$ in the $\bxi$ coordinate system~\cite{Hutton2003}.

The introduction of the volume coordinate system has two main purposes: first, it allows a simple description of polynomial basis functions of any order, reducing to $\bvarphi(\bxi)=\bxi$ for first order FEM~\cite{Beltrachini2015a}; and second, it simplifies the solution of integrals of the form
\begin{equation}\label{integ_vol_coord}
\int_{\mathcal{T}} \xi_1^i \xi_2^j \xi_3^k \xi_4^l d\bxi=6V\frac{i!j!k!l!}{(i+j+k+l+3)!}.
\end{equation}
After transforming~\eqref{eq_ke} to the volume coordinate system, and noting that the gradient of any scalar function $h(\br)$ can be expressed in the volumetric coordinate system as $\nabla h(\br) = (6V)^{-1} \bLambda \nabla_{\xi} h(\bxi)$~\cite{Beltrachini2015a}, we obtain
\begin{equation} \label{Se}
\bK^e=\frac{1}{6V}\int_{\mathcal{T}_{n}} \nabla_{\xi}\bvarphi^T(\bxi)\ \bLambda^T\overline{\bsigma}\bLambda\ \nabla_{\xi}\bvarphi(\bxi) d\bxi,
\end{equation}
where $\nabla_{\xi}\bvarphi(\bxi)=\big[ \nabla_{\xi}\varphi_1(\bxi),\dots,\nabla_{\xi}\varphi_{N}(\bxi) \big]$ and $\nabla_\xi$ is the gradient operator in the volume coordinate space. To further simplify~\eqref{Se}, we apply the $vec$ operator to $\bK^e$, i.e. the linear operator that transforms the element matrix in a column vector obtained by stacking the columns of $\bK^e$ on top of one another. Using the identity $\vect\left( \bA \bB \bC \right)=\big( \bC^T \otimes \bA \big) \vect(\bB)$~\cite{Horn1991}, we get
\begin{equation*}
\vect \left( \bK^e \right)= \frac{1}{6V} \bQ^T \left(\bLambda\otimes \bLambda \right)^T \ \vect(\overline{\bsigma}),
\end{equation*}
where 
\begin{equation*}
\bQ=\int_{\mathcal{T}_{n}} \big(\nabla_{\xi}\bvarphi(\bxi)\otimes \nabla_{\xi}\bvarphi(\bxi)\big) d\bxi,
\end{equation*}
is a matrix with constant coefficients, and $\otimes$ is the Kronecker product. Once the basis functions are chosen, $\bQ$ is easily computed using~\eqref{integ_vol_coord}. In the particular case of using first order basis functions, $\nabla_{\xi}\bvarphi(\bxi)=\bI_4$, and consequently $\bQ=\bI_{16}/6$. Therefore, the element matrix considering first order basis functions turns out to be
\begin{equation}\label{Ke_parametric}
\bK^e=\frac{1}{36 V} \bLambda^T\overline{\bsigma}\bLambda.
\end{equation}

\section{Source vector}\label{sec:4}

The non-linearities introduced by $f(\br)$ impose new difficulties in the computation of the element source vectors. To overcome them, we employ Gauss theorems defined on local reference frames to reduce their dimensionality to 1D. In this Section, we first review and extend such framework~\cite{Wilton1984,Graglia1993}, for then utilising it in the computation of~\eqref{eq_bv} and~\eqref{eq_bs}.

\subsection{Local coordinate system}\label{sec:4loc}

\begin{figure}[th]
\centering
\includegraphics[width=120mm]{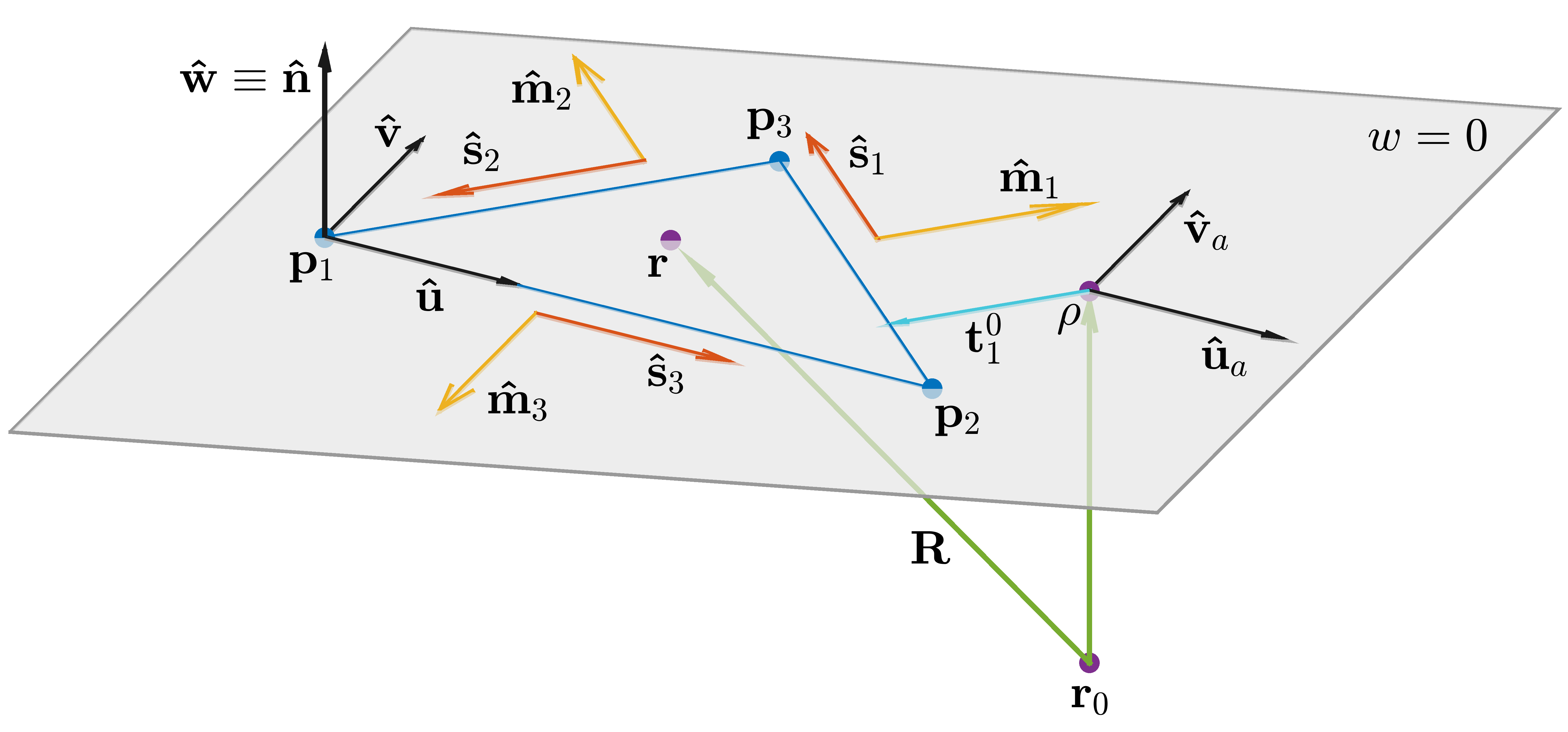}
\caption{Schematic representation of the local frames $\mathcal{L}$ and $\mathcal{L}_a$ (see Section~\ref{sec:4loc}).}
\label{fig_0}
\end{figure}

Let $\mathcal{G}\equiv [x,y,z]$ be the global Cartesian frame in which our model is inscribed (see Fig.~\ref{fig_0}). We consider a triangle $T$ formed by the nodes $\bp_1$, $\bp_2$, and $\bp_3$, such that $\partial T_i$ represents its $ith$ side, i.e. opposite to $\bp_i$ ($i=1,2,3$). For this triangle, we define the vectors
\[
\bs_i=\bp_{i-1}-\bp_{i+1}, \ i=1,2,3,
\]
with $\bp_{4}\equiv \bp_1$ and $\bp_{-1}\equiv \bp_3$. These vectors allow to define a local frame $\mathcal{L}\equiv [u,v,w]$, with origin in $\bp_1$, such that $T$ belongs to the plane $w=0$. The unitary vectors for the local frame are chosen as
\[
\hat{\bu}=\hat{\bs}_3,\ \ \hat{\bv}=\hat{\bw} \times \hat{\bu},\ \ \hat{\bw}=\hat{\bs}_1 \times \hat{\bs}_{2},
\]
where $\hat{\bs}_i=\bs_i/s_i$, with $s_i=|\bs_i|$, $i=1,2,3$. In this new frame, any point $\br\in T$ will have coordinates $[u,v,0]_\mathcal{L}\equiv [u,v]$. In particular, 
\begin{equation}\label{eq:points}
\bp_1=[0,0],\ \bp_2=[s_3,0], \ \bp_3=[u_3,v_3],
\end{equation}
with $u_3=-\bs_2\cdot \hat{\bu}$ and $v_3=-\bs_2\cdot\hat{\bv}$.

We are interested in solving integrals involving the vector $\bR$ and/or its magnitude $R$. In the local frame, the source position can be written as $\br_0=\brho - w_0 \hat{\bw}$, where $\brho=[u_0,v_0,0]_{\mathcal{L}}$ is the projection of $\br_0$ onto the plane $w=0$. To perform such integrals, it is convenient to translate $\mathcal{L}$ so that the origin is $\brho$. This new frame, called $\mathcal{L}_a\equiv [u_a,v_a,w]$, has the same unitary vectors as $\mathcal{L}$. Noting that $u=u_a+u_0$ and $v=v_a+v_0$, it can be easily seen that 
\begin{gather*}
\br=[u,v,0]_{\mathcal{L}}=[u_a,v_a,0]_{\mathcal{L}_a},\\
\br_0=[u_0,v_0,-w_0]_{\mathcal{L}}=[0,0,-w_0]_{\mathcal{L}_a},\\
\bR=[u-u_0,v-v_0,w_0]_{\mathcal{L}}=[u_a,v_a,w_0]_{\mathcal{L}_a}.
\end{gather*}

The use of $\mathcal{L}_a$ facilitates the parameterisation of any point in $\partial T$, which will be required for calculating the line integrals after the application of Gauss theorems. Let $\bt_i^0$ be the shortest vector going from $\brho$ to $\partial T_i$. These vectors can be written as 
$\bt_i^0=t_i^0 \hat{\bm}_i$, with 
\[
\hat{\bm}_i=\hat{\bs}_i\times \hat{\bw}, \ i=1,2,3.
\]
The coefficients $t_i^0$ can be obtained from the intersection between the lines $(\brho,\hat{\bm}_i)$ and $(\bp_{i-1},\hat{\bs}_i)$, resulting in~\cite{Graglia1993}
\begin{align*}
&t_1^0=\left[ v_0(u_3-s_3)+v_3 (s_3-u_0) \right]/s_1,\\
&t_2^0=\left( u_0 v_3 - v_0 u_3 \right)/s_2,\\
&t_3^0=v_0.
\end{align*}
Then, every point $\bs\in\partial T_i$ can be parameterised as $\bs(s)=\bt_i^0+s \hat{\bs}_i$, where $s\in [s_i^-,s_i^+ ]$ satisfies $\bs(s_i^-)=\bp_{i-1}$ and $\bs(s_i^+)=\bp_{i+1}$. After some algebra, we find~\cite{Graglia1993}
\begin{align*}
&s_1^-=-\left[ (s_3-u_0)(s_3-u_3)+v_0 v_3 \right]/s_1,\\
&s_1^+=\left[ (u_3-u_0)(u_3-s_3)+ v_3 (v_3-v_0) \right]/s_1,\\
&s_2^-=-\left[ s_3(s_3-u_0)+v_3 (v_3-v_0) \right]/s_2,\\
&s_2^+=\left( u_0 u_3 + v_0 v_3 \right)/s_2,\\
&s_3^-=-u_0,\ \ \ s_3^+=s_3-u_0.
\end{align*}

In the following, $R_i^0=\sqrt{(t_i^0)^2+w_0^2}$ is the distance between $\br_0$ and $\partial T_i$, and $R_i^\pm=\sqrt{(s_i^\pm)^2+(t_i^0)^2+w_0^2}$ is the distance between $\br_0$ and $\bs(s_i^\pm)$, $i=1,2,3$. Also, to express results concisely, we define the following functions (some of them introduced in~\cite{Wilton1984,Graglia1993})
\begin{gather*}
f_{2i}=\ln\left( \frac{R_i^+ + s_i^+}{R_i^- + s_i^-} \right),\\
R^s_i= \frac{1}{\left( R_i^0\right)^2} \left( \frac{s_i^+}{R_i^+}-\frac{s_i^-}{R_i^-} \right),\\
\beta_i=\tan^{-1} \frac{t_i^0 s_i^+}{\left(R_i^0\right)^2 + |w_0| R_i^+}  - \tan^{-1} \frac{t_i^0 s_i^-}{\left(R_i^0\right)^2 + |w_0| R_i^-},
\end{gather*}
with the arctangent function evaluated on its principal branch.

\subsection{Gauss theorems}

We are interested in the application of Gauss theorems for reducing the dimensionality of 2D and 3D integrals to 1D, where analytical expressions can be easily found. Let $S$ be a planar surface, and $\hat{\bu}_m$ the unitary vector normal to its boundary $\partial S$ and belonging to the same plane (pointing outwards). Then, for any continuously differentiable functions $g:S\rightarrow \mathbb{R}$ and $\bg:S\rightarrow \mathbb{R}^3$, the surface Gauss theorems state that~\cite{VanBladel2007}
\begin{equation}\label{eq:gauss_2d}
\int_S \nabla_S \cdot \bg \ dS = \int_{\partial S} \bg \cdot \hat{\bu}_m \ ds, \ \ \int_S \nabla_S g \ dS = \int_{\partial S} g \hat{\bu}_m \ ds,
\end{equation}
where $\nabla_S=\hat{\bu} \partial /\partial u  + \hat{\bv} \partial /\partial v $ is the {\it surface vector differential operator}. Similarly, let $V$ be a volume and $\hat{\bu}_n$ the unitary vector normal to $\partial V$, the boundary of $V$ (pointing outwards). Then, for any functions $g:V\rightarrow \mathbb{R}$ and $\bg:V\rightarrow \mathbb{R}^3$, the volume Gauss theorems are given by~\cite{VanBladel2007}
\begin{equation}\label{eq:gauss_3d}
\int_V \nabla \cdot \bg \ dV = \int_{\partial V} \bg \cdot \hat{\bu}_n \ dS, \ \ \int_V \nabla g \ dS = \int_{\partial V} g \hat{\bu}_n \ dS.   
\end{equation}

In the present study, since we are assuming tetrahedral meshes, the application of~\eqref{eq:gauss_3d} to volume integrals will result in transformed integrals over their triangular faces. Therefore, special emphasis should be put in expressing the integrands in a proper differential form, so that~\eqref{eq:gauss_2d} can be utilised. To this end, the following proposition (shown in Appendix~1) will be of practical use.

\begin{proposition}\label{prop_1} Let $h:S\rightarrow \mathbb{R}$ be a function defined for any vector $\bt$ belonging to the planar surface~$S$. Then, the following relation holds
\[
\nabla_S\cdot \left( h(R) \frac{\bt}{t^2}\right)=\frac{\bt}{t^2}\cdot\nabla_Sh(R).
\]\vspace{0.1mm}
\end{proposition}

Proposition~\ref{prop_1} allows to express integrands in the form required to apply Gauss theorems on a planar surface. In particular, if we choose $h(R)=R^{-1}$ and $h(R)=R^{-3}$, we get 
\begin{align}
&\frac{1}{R^3}=-\nabla_S\cdot \left(\frac{1}{R} \frac{\bt}{t^2}\right),\label{eq:R3}\\
&\frac{3}{R^5}=-\nabla_S\cdot \left(\frac{1}{R^3} \frac{\bt}{t^2}\right),\label{eq:R5}
\end{align}
respectively. In the special case of choosing $\bt$ as $\hat{\bu}$ or $\hat{\bv}$, these expressions turn out to be
\begin{align}
&\begin{Bmatrix}
u_a\\
v_a
\end{Bmatrix}
\frac{1}{R^3}
=-\nabla_S \cdot \left(
\begin{Bmatrix}
\hat{\bu}\\
\hat{\bv}
\end{Bmatrix}
 \frac{1}{R} \right), \label{eq:R3a}\\
&
\begin{Bmatrix}
u_a\\
v_a
\end{Bmatrix}
\frac{3}{R^5}=-\nabla_S \cdot \left(
\begin{Bmatrix}
\hat{\bu}\\
\hat{\bv}
\end{Bmatrix}
\frac{1}{R^3}\right), \label{eq:R5a}  
\end{align}
where the braces indicate that the relation is valid for either the upper or lower elements in them. Finally, utilising simple differentiation rules, we obtain
\begin{equation}\label{eq:Rb}
\begin{Bmatrix}
u_a\\
v_a
\end{Bmatrix} 
\nabla_S \frac{1}{R^{3}} = \nabla_S \left( 
\begin{Bmatrix}
u_a\\
v_a
\end{Bmatrix}
\frac{1}{R^3} \right)  - 
\begin{Bmatrix}
\hat{\bu}\\
\hat{\bv}
\end{Bmatrix}
\frac{1}{R^3}.
\end{equation}

\subsection{Surface integral}\label{sec:surf}

Now, we take advantage of the results from the previous subsections to compute the surface element vector~$\bb^s$. In this case, the linear interpolation vector is $\bvarphi(\br)=[\varphi_1(\br),\varphi_2(\br),\varphi_3(\br)]^T$, with $\varphi_i (\bp_j)=\delta_{ij}$, for every $i,j=1,2,3$. In the $\mathcal{L}$ frame, $\bvarphi$ takes the form~\cite{Zienkiewicz2013}
\begin{equation}\label{eq:phi}
\bvarphi(u,v)
=
\begin{bmatrix}
1 & -s_3^{-1} & \left(u_3 s_3^{-1}-1\right)v_3^{-1} \\
0 & s_3^{-1} & -u_3 s_3^{-1} v_3^{-1} \\
0 & 0 & v_3^{-1}
\end{bmatrix}
\begin{bmatrix}
1 \\
u \\
v 
\end{bmatrix},
\end{equation}
which clearly complies the interpolation conditions~\eqref{eq:points}. Replacing~\eqref{eq:phi} in~\eqref{eq_bs} leads to
\begin{equation}\label{eq:su1}
4 \pi\bb^s=\bvarphi(u_0,v_0) I_0 + 
\begin{bmatrix}
-s_3^{-1} & \left(u_3 s_3^{-1}-1\right)v_3^{-1} \\
s_3^{-1} & -u_3 s_3^{-1} v_3^{-1} \\
0 & v_3^{-1}
\end{bmatrix}
\begin{bmatrix}
I_{u}^a \\
I_{v}^a
\end{bmatrix},
\end{equation}
where 
\begin{align}
&I_0= \int_{T} \hat{\bn} \cdot \nabla f(\br) d\br, \label{eq:su2}\\
&\begin{Bmatrix}
I_u^a \\
I_v^a
\end{Bmatrix}=\int_{T} 
\begin{Bmatrix}
u_a \\
v_a
\end{Bmatrix} \hat{\bn} \cdot \nabla f(\br) d\br. \label{eq:su3}
\end{align}
Then, we need to find $I_0$, $I_u^a$, and $I_v^a$ for any arbitrary dipolar source. To do so, the following proposition (shown in Appendix~2) is of great help.

\begin{proposition}\label{prop_2} For any triangle $T$ and dipolar source located in $\br_0$ with moment $\bq$, the following relation holds
\[
\hat{\bn} \cdot \nabla \left( \bq \cdot \frac{\bR}{R^3} \right)= \bq \cdot \nabla \left( \hat{\bn} \cdot \frac{\bR}{R^3} \right).
\] \vspace{.01cm}
\end{proposition}

Proposition~\ref{prop_2} is used to rewrite~\eqref{eq:su2}~--~\eqref{eq:su3} as
\begin{align}
&I_0= \bq \cdot \int_{T} \nabla \left( \hat{\bn} \cdot \frac{\bR}{R^3} \right) d\br,\label{eq:su4}\\
&\begin{Bmatrix}
I_u^a \\
I_v^a
\end{Bmatrix}
=\bq \cdot\int_{T} 
\begin{Bmatrix}
u_a \\
v_a
\end{Bmatrix}  
\nabla \left( \hat{\bn} \cdot \frac{\bR}{R^3} \right) d\br, \label{eq:su5}
\end{align}
which allows to exploit the fact that, in~$\mathcal{L}_a$, $\hat{\bn}=\hat{\bw}$, resulting in $\hat{\bn}\cdot \bR = w_0$. 

To solve~\eqref{eq:su4} in the $\mathcal{L}_a$ frame, we first note that, for any function $g$ defined in $\mathcal{L}_a$, $\nabla g = \nabla_S g + \partial g / \partial w_0 \hat{\bw}$. Applying this property to the integrand of~\eqref{eq:su4}, we obtain   
\begin{equation}\label{eq:I_dec}
\nabla \frac{w_0}{R^3}=\nabla_S \frac{w_0}{R^3}+\frac{\hat{\bw}}{R^3}-3\frac{w_0^2}{R^5}\hat{\bw}.
\end{equation}
Eq.~\eqref{eq:I_dec} allows us to split $I_0$ in three terms, which are treated separately. The first term can be easily solved by applying the second identity from~\eqref{eq:gauss_2d}, resulting in
\begin{align}
\int_T \nabla_S \frac{w_0}{R^3} dT &= \int_{\partial T} \frac{w_0}{R^3} \hat{\bm} ds = w_0 \sum_{i=1}^3 \hat{\bm}_i \int_{\partial T_i} \frac{ds_i}{R^3} \nonumber \\
&=w_0 \sum_{i=1}^3 \hat{\bm}_i \int_{s_i^-}^{s_i^+} \frac{ds_i}{\left( w_0^2 + (t_i^0)^2 + s_i^2  \right)^{3/2}} \nonumber \\ &= w_0 \sum_{i=1}^3 \hat{\bm}_i R^s_i. \label{eq:I0_1}
\end{align}
The computation of the second and third terms is done using~\eqref{eq:R3} and~\eqref{eq:R5}, after which the first relation in~\eqref{eq:gauss_2d} is applied. However, it should be noted that the utilisation of~\eqref{eq:R3} and~\eqref{eq:R5} introduce a singularity in $\brho$ (since~$t=0$ for~$\brho$ in~$\mathcal{L}_a$), which prevents the application of Gauss theorems straightforwardly. This issue can be solved by splitting the domain in two parts, one being $T_\epsilon$, a circle with radius $\epsilon$ centred in $\brho$, and the other its complement in $T$, i.e. $T-T_\epsilon$. Then, the corresponding integrals are computed separately, after which the limit for $\epsilon \rightarrow 0$ is performed. This procedure was followed by Graglia~\cite{Graglia1993}, who obtained the following expression for the second term 
\begin{equation}
\int_T \frac{dT}{R^3}=\frac{\bbeta}{|w_0|},\label{eq:I0_2}
\end{equation}
where $\bbeta=\beta_1+\beta_2+\beta_3$. In Appendix~3  we show that the integral corresponding to third term is given by
\begin{equation}
\int_T \frac{dT}{R^5}=\frac{1}{3}\left( \frac{\bbeta}{|w_0|^3} + \sum_{i=1}^3 \frac{t_i^0}{w_0^2}R_i^s\right).\label{eq:I0_3}
\end{equation}
Finally, the proper combination of~\eqref{eq:I0_1}~--~\eqref{eq:I0_3} yields
\begin{equation}\label{eq:IR5}
I_0= \bq \cdot \sum_{i=1}^3 R^s_i \left( w_0  \hat{\bm_i}  -  t_i^0  \hat{\bw}\right).
\end{equation}

The integrals $I^a_u$ and $I^a_v$ are obtained following a similar procedure. In Appendix~4 we show that
\begin{align}
&\begin{Bmatrix}
I_u^a\\
I_v^a
\end{Bmatrix}
=
\bq \cdot \Bigg[ w_0 \sum_{i=1}^3 \hat{\bm}_i 
\begin{Bmatrix}
\hat{\bu}\\
\hat{\bv}
\end{Bmatrix}
\cdot \left(\hat{\bs}_i R_i^d + \bt_i^0 R_i^s \right) \nonumber \\ 
&\hspace{3mm}-
\begin{Bmatrix}
\hat{\bu}\\
\hat{\bv}
\end{Bmatrix}
\sign(w_0)\bbeta 
+ \hat{\bw} \left(
\begin{Bmatrix}
\hat{\bu}\\
\hat{\bv}
\end{Bmatrix}
\cdot \sum_{i=1}^3 \hat{\bm}_i \left( R_i^s w_0^2 - f_{2i} \right)
\right) \Bigg], \label{eq:Iuava}
\end{align}
where $R_i^d=(R_i^-)^{-1} - (R_i^+)^{-1}$ and $\sgn$ is the sign function.

\subsection{Volume integral}

The computation of $\bb^v$ turns out to be much simpler than the corresponding to the surface element vector. Utilising~\eqref{eq:gauss_3d} and the fact that $\nabla \bvarphi=(6V)^{-1} \bLambda$ for first order basis functions (see Section~\ref{sec:3}), \eqref{eq_bv} reduces to
\[
\bb^v=\frac{\bLambda^T \overline{\bsigma}^c}{24 V \pi \sigma^\infty} \sum_{i=1}^4 \hat{\bn}_i \int_{T^i} \nabla f d\br,
\] 
where $T^i$ is the $ith$ triangular face of the tetrahedron under consideration. Noting that, for a dipolar source, $f(\br)=-\bq \cdot \nabla R^{-1}$, we get 
\begin{align}
\int_{T^i} \nabla f d\br&=-\bq\cdot \int_{T^i} \nabla \frac{1}{R} d\br \nonumber\\
&= \bq\cdot \left(-\hat{\bw}^i\sgn(w_0^i) \bbeta^i + \sum_{l=1}^3 \hat{\bm}_l^i f_{2l}^i \right),
\end{align}
where the last equality follows from~\cite[eq. (34)]{Graglia1993} (and noting that, in this paper, $\bw_0$ is opposite to its analogous in~\cite{Graglia1993}) and the superscripts refer to the triangular face under consideration.

\section{Experiments}\label{sec:6}

\subsection{Local errors}

We first studied the performance of the AS method for computing the element vectors $\bb^s$ and~$\bb^v$. This was done by calculating the analytical solution for a given element and source of electrical activity, and comparing the results with those obtained with the Gauss-Jacobi quadrature scheme, known as the {\it full subtraction} (FS) approach~\cite{Drechsler2009,Beltrachini2019b}. More specifically, we considered a source in $\br_0=[d,0,0]$ with moment $\bq=[10,0,0]nAm$ and calculated the element vectors for equilateral simplices located in the origin with side length $a$. In the case of $\bb^s$, we used a triangle with vertices $\bp_1=[0,-a/2,-a\sqrt{3}/6]$, $\bp_2=[0,a/2,-a\sqrt{3}/6]$, and $\bp_3=[0,0,a\sqrt{3}/3]$; in the case of $\bb^v$, we used a tetrahedron with nodes $\bp_1=[0,0,a\sqrt{3}/3]$, $\bp_2=[0,a/2,-a\sqrt{3}/6]$, $\bp_3=[0,-a/2,-a\sqrt{3}/6]$, and $\bp_4=[-a\sqrt{3}/3,0,0]$. Then, we computed the relative error (RE) between element vectors, defined as $RE_e=\|\bb_n - \bb_a\| / \| \bb_a \|$, with $\bb_n$ and $\bb_a$ being the numerical and analytical solutions, respectively. We repeated such experiments for different combinations of $d$ and $a$.

\subsection{Spherical head model}

In a second experiment, we evaluated the improvements of the AS approach over numerical quadrature schemes in the solution of the EEG-FP. To do so, we utilised a multi-layered spherical head model with anisotropic electrical conductivity field, for which analytical solutions are available~\cite{deMunck1993,Beltrachini2019b}. We considered four compartments representing the scalp, skull, cerebrospinal fluid (CSF), and brain, each with outer radius equal to 0.092m, 0.086m, 0.08m, and 0.078m, respectively. The electrical conductivities were considered isotropic for the scalp, CSF, and brain compartments, and set to 0.33 S/m, 1.79 S/m, and 0.33 S/m, respectively. The skull was modelled as an anisotropic layer with radial/tangential conductivities of 0.0093/0.015 S/m. The electric conductivity values were extracted from the relevant literature~\cite{Baumann1997,Ramon2006,Dannhauer2011,McCann2019}.

We discretised the spherical model using the ISO2Mesh toolbox~\cite{Fang2009}. Meshes were built to achieve a maximum {\it radius-edge factor} of~1.2. We considered a coarser mesh resolution in the brain layer since the volume integral in~\eqref{eq_l} vanishes in there for sharing the same electrical conductivity as the source neighbourhood. Six discretisations with 39k, 119k, 281k, 440k, 640k, and 938k nodes were obtained. A total of 162 electrodes uniformly placed on the scalp surface were utilised~\cite{Koay2011}. A visual representation of the models is available in the Supplementary material.

We simulated radially- and tangentially-oriented sources in~100 different locations at a distance~$r_0$ from the next conductivity jump, with~$r_0$ taking values in the range $0.125-1.5mm$. For each of them, we computed the relative error (RE), relative difference measure (RDM), and magnification error (MAG), defined as
\begin{align*}
&RE=\frac{\|\bu_n - \bu_a\|}{ \| \bu_a \|},\\
&RDM=\left\| \frac{\bu_n}{ \|\bu_n\|} - \frac{\bu_a} {\| \bu_a \|}\right\| ,\\
&MAG=\left|1- \frac{\| \bu_n \| }{ \|\bu_a \|} \right|,
\end{align*}
with $\bu_n$ and $\bu_a$ being the numerical and analytical solutions of the EEG-FP, respectively. The RE is a measure of the overall difference between the analytical and numerical solutions of the EEG-FP, whereas the RDM and MAG account for topographical and magnitude differences between them, respectively. Error measures were computed for the solutions obtained with the AS and FS approaches, the latter with $n=2$ and $n=4$. 

In addition to standard error measures, we computed the RE between the numerical approximations provided by the AS and FS methods, 
\begin{equation*}
RE_s=\frac{\|\bu_{as} - \bu_{fs}\|}{ \| \bu_{as} \|},
\end{equation*}
where $\bu_{as}$ and $\bu_{fs}$ are the numerical solutions computed with the AS and FS approaches, respectively. This allowed to evaluate the errors introduced by the utilisation of a numerical quadrature scheme, which are avoided using the AS approach. We calculated the~$RE_s$ for all mesh discretisations considering the aforementioned sources.

\subsection{Illustration in a real model}\label{sec:real_1}

Finally, we illustrated the advantages of the AS method over the FS approach in a personalised head model. We constructed a realistic head representation based on the Colin~27 high resolution MRI segmentation of the Montreal Neurological Institute~\cite{Aubert2006}. A mesh with 3.8M tetrahedral elements (600k nodes) was created using the ISO2Mesh software as described before. This included the scalp, skull, CSF, grey matter (GM), and white matter (WM) compartments, each of them assumed isotropic and with electrical conductivities equal to 0.33~S/m, 0.004~S/m, 1.79~S/m, 0.45~S/m, and 0.13~S/m, respectively~\cite{McCann2019}. A total of 162 electrodes were placed on the scalp according to the ABC standard. A visual representation of the model is available in the Supplementary material.

We computed the leadfield matrix for sources on a surface located at the mid-point between the GM/CSF and WM/GM layers, as provided by FreeSurfer~\cite{Dale1999}. Since the cortical thickness is in the range 1-5mm~\cite{Fischl2000}, extra mesh refinement would be needed in thin cortical regions to avoid the FS method to become unstable. We calculated the leadfield matrices with the AS and FS ($n=2$) techniques and obtained the $RE_s$. As in the spherical case, this allowed to evaluate the errors introduced in the solution of the EEG-FP due to the use of numerical quadrature.

\section{Results}\label{sec:7}

\subsection{Local errors}

\begin{figure}[htbp]
\centering
\includegraphics[width=11cm]{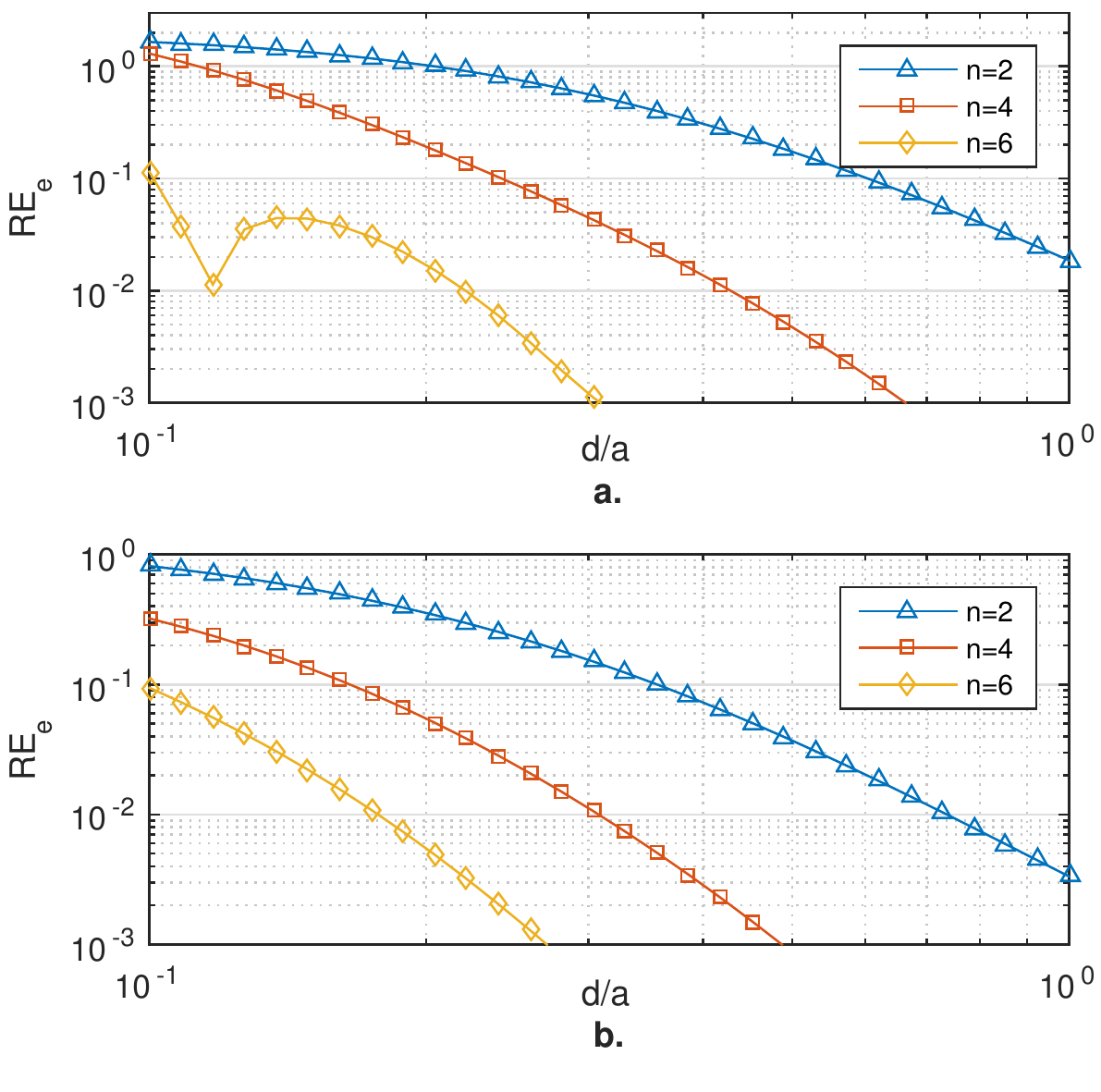}
\caption{RE between the analytical and numerical element source vectors as a function of the normalised distance between the source and the element. Results are presented for the surface (a.) and volume (b.) element vectors and numerical integration orders $n$ equal to 2, 4, and 6 (with different markers).}
\label{fig_1}
\end{figure}

The REs between the element vectors obtained with the analytical and numerical approaches are shown in Fig.~\ref{fig_1}. Figs.~\ref{fig_1}a. and b. display the $RE_e$ for the surface and volume element source vectors, respectively, as a function of $d/a$ and for integration orders 2, 4, and 6. In addition to providing more accurate solutions, the analytical method allowed to reduce the computational demands on their calculation. In the case of the surface element vector, the analytical method resulted 2, 6, and 12 times faster compared to the numerical technique with~$n$ equal to 2, 4, and 6, respectively (in average over 500 iterations). Regarding the volume element vector, the speed-up of the analytical method was 1.4, 4.2, and 11.1 for the same integration orders and number of iterations. It is worth mentioning that the same results were obtained for a given~$d/a$ independently of the values chosen for~$a$ and~$d$, highlighting the importance of such ratio in the error introduced by the FS methodology. 

\subsection{Spherical head model}

\begin{figure*}[htbp]
\centering
\includegraphics[width=17cm]{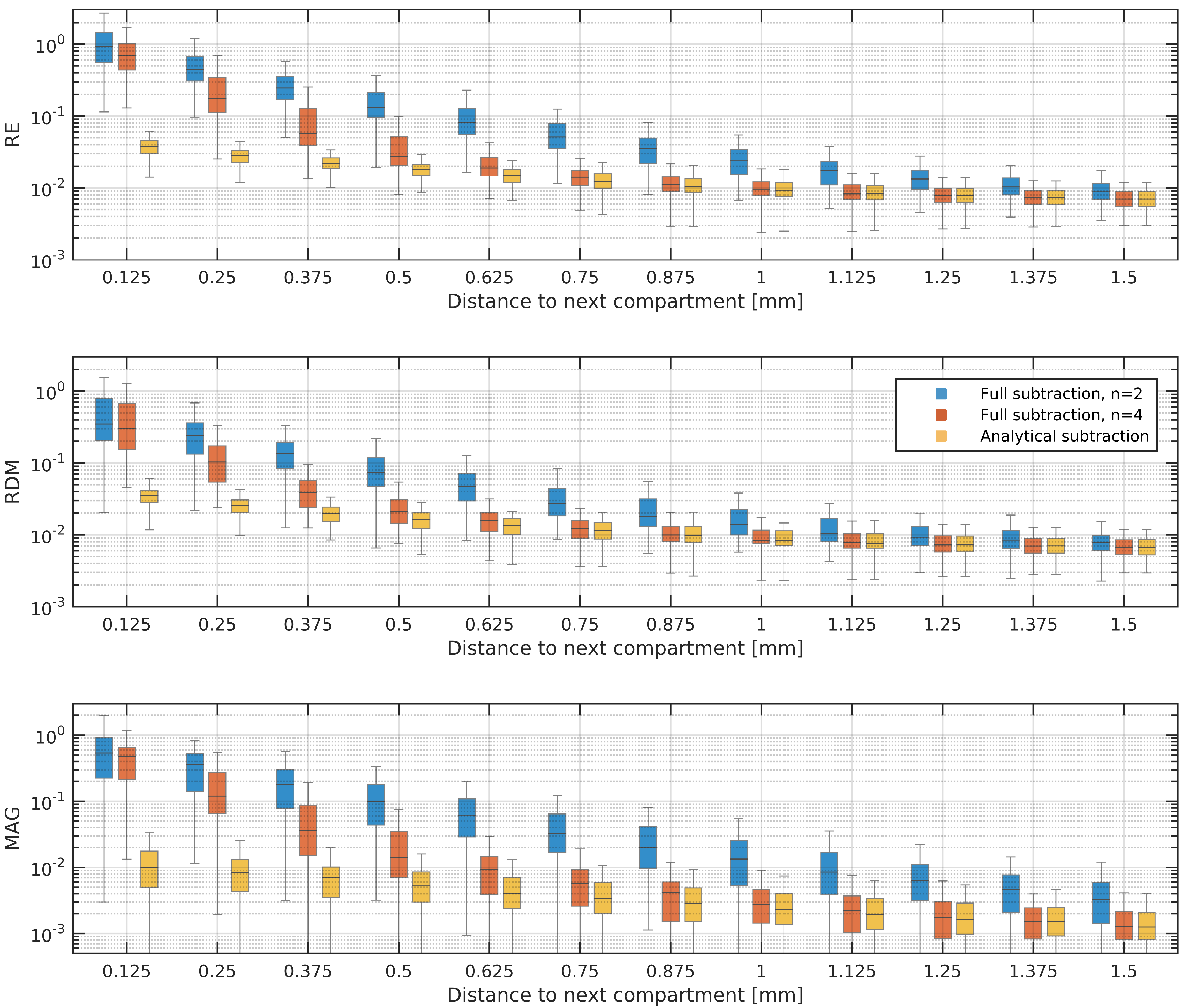}
\caption{Error measures (RE: top; RDM: centre; MAG: bottom) for the numerical solutions of the EEG-FP as a function of the distance to the next compartment considering tangentially-oriented sources and the model with 256k nodes. Results are presented for the AS and FS (with $n=2$ and $n=4$) methods (with different colours).}
\label{fig_2}
\end{figure*}

In Fig.~\ref{fig_2} we show the RE, RDM, and MAG as a function of the distance to the closest conductivity change considering the model with 256k nodes and the AS and FS methods, the latter with $n=2$ and $n=4$. Results are displayed for tangentially-oriented and highly eccentric sources for which the numerical method was shown to introduce noticeable errors. Error measures for radially-oriented sources and other models are presented in the Supplementary material. Regarding the computational cost of the methods, the AS approach resulted, in average, 1.04 and 2.11 times faster than using the FS technique with $n$ equal to 2 and 4, respectively.

\begin{figure*}[htbp]
\centering
\includegraphics[trim=2cm 0cm 0.5cm 0cm,clip,width=16cm]{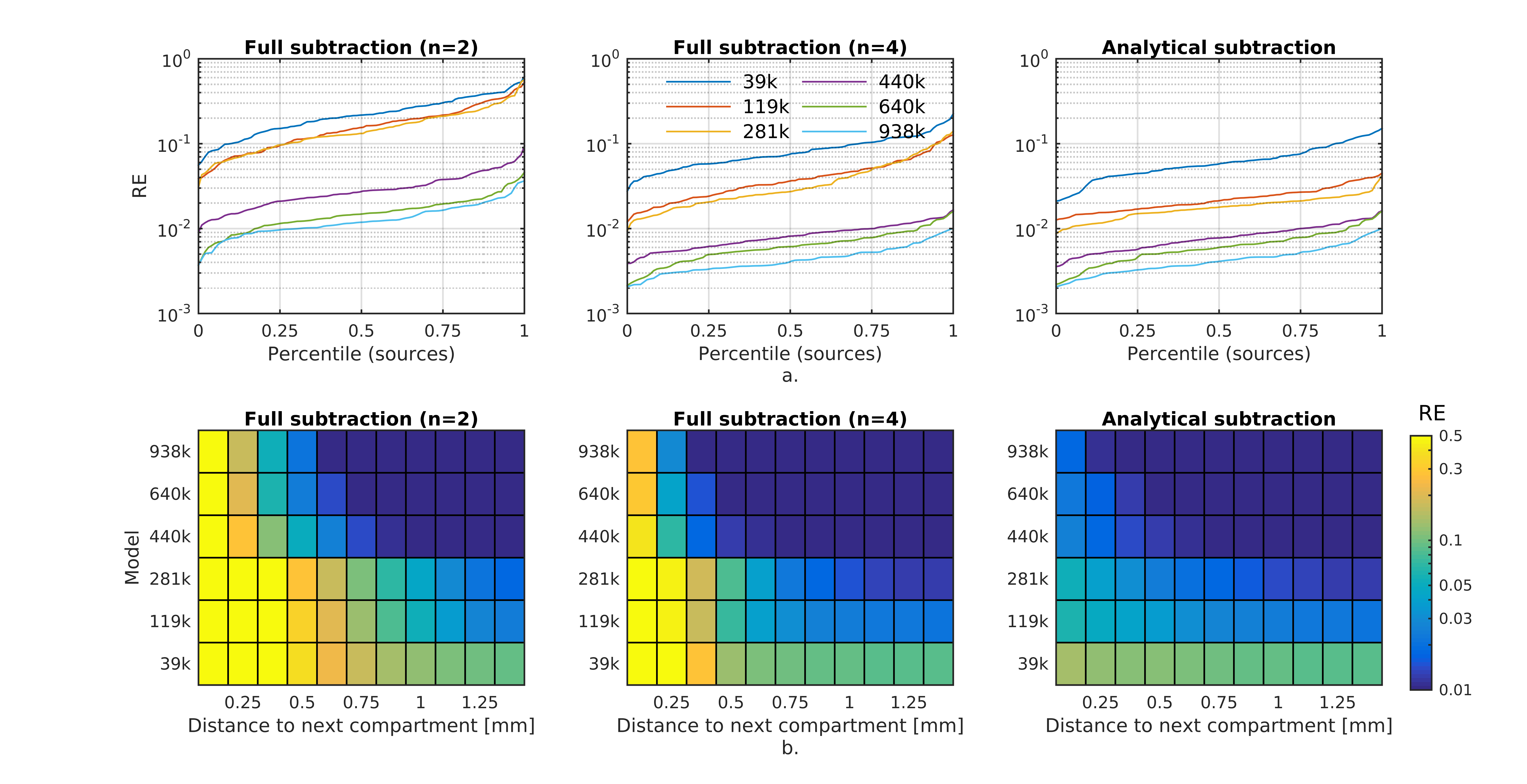}
\caption{RE for different models and source locations. a. Percentile curves of the RE considering 100 sources located at 0.5mm to the next compartment using all models (with different colours) and the FS approach with $n=2$ (left), $n=4$ (centre), and the AS method (right). b. 90th percentile of the RE as a function of the distance to the next compartment and the number of nodes. Results considering the FS (with $n=2$ and $n=4$) and AS methods are displayed as in a.}
\label{fig_3}
\end{figure*}

Fig.~\ref{fig_3} displays the percentiles of the RE obtained for different source locations and model discretisations. In Fig.~\ref{fig_3}a. we show the percentile curves of the RE considering 100 tangentially-oriented sources located at 0.5mm to the next compartment and all models. As expected, the errors obtained with the AS approach are always lower or equal to those obtained with the FS method regardless of the integration order and head model. Such difference becomes more apparent as the source is closer to the next compartment and the model is less refined. Fig.~\ref{fig_3}b. presents the 90th percentile of the RE as a function of the distance to the next compartment and the number of first order nodes for both AS and FS ($n=2$ and $n=4$) techniques. It is evident that the FS approach provides solutions with unacceptable REs as the source becomes closer to the closest conductivity jump. Although the increase of the integration order allows to reduce the distance for which the method becomes unstable, this issue cannot be completely avoided by using a numerical quadrature scheme. This is not the case of the AS method, which provides stable solutions for arbitrarily eccentric sources.   

\begin{figure*}[htbp]
\centering
\includegraphics[width=11cm]{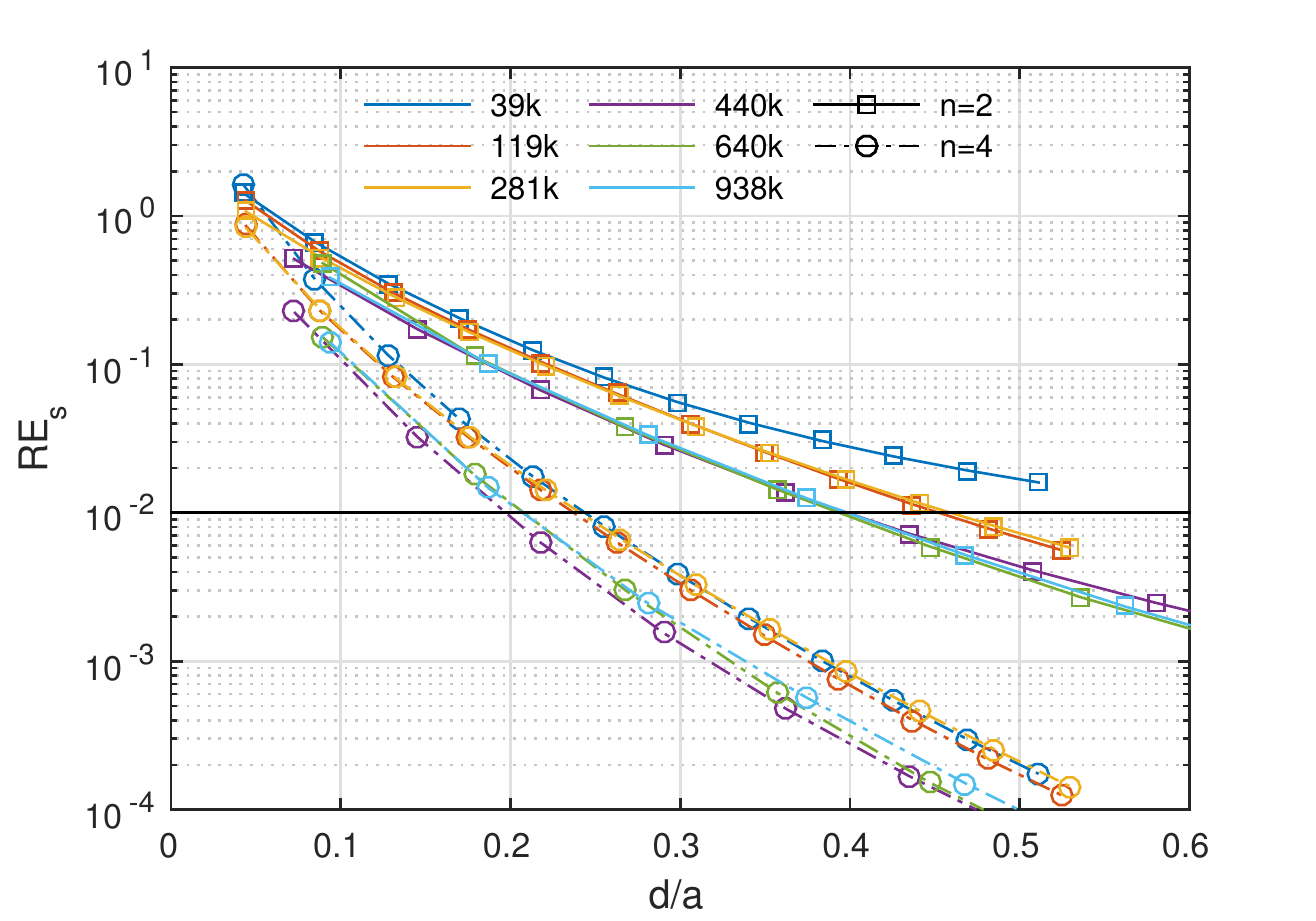}
\caption{Mean RE between the numerical approximations obtained with the AS and FS approaches as a function of $d/a$. Results are presented for all mesh discretisations (with different colours) and numerical integration orders $n=2$ and~$n=4$ (with different markers). The black horizontal line indicates the $1\%$ threshold (see text).}
\label{fig_4}
\end{figure*}

The mean RE between the AS and FS approximations is presented in Fig.~\ref{fig_4}. There, we plot the $RE_s$ as a function of the ratio $d/a$, considering all model discretisations and integration orders $n=2$ and~$n=4$. We considered $a$ as the mean element side-length belonging to the CSF compartment (i.e. those belonging to the closest compartment with different electrical conductivity). It can be noted that the approximations obtained with the FS method differ from those obtained with the AS approach consistently as a function of $d/a$ regardless of the model discretisation used. In the case of adopting the FS method with $n=2$, errors due to the use of a numerical integration scheme start becoming noticeable (i.e. larger than $1\%$) for $d/a<0.5$. This threshold reduces to $d/a<0.25$ if the integration order is increased to $n=4$.  

\subsection{Real model}

\begin{figure*}[htbp]
\centering
\includegraphics[trim=1cm 0cm 1.2cm 0cm,width=16cm]{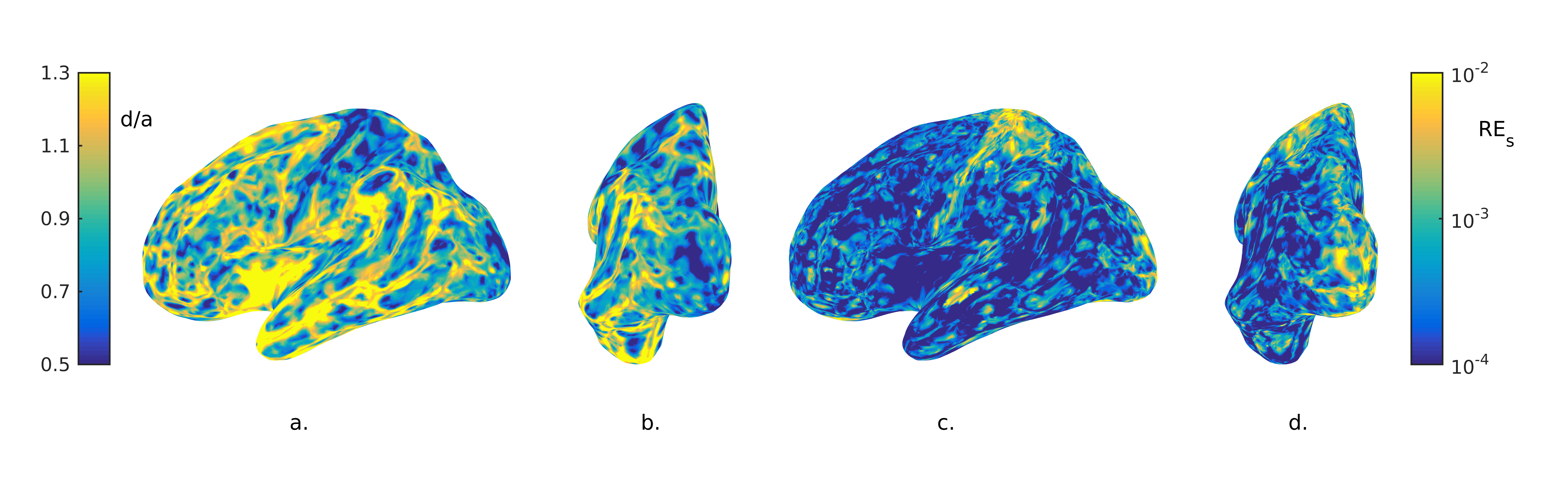}
\caption{a--b. Ratio $d/a$ for the mesh discretisation and sources described in Section~\ref{sec:real_1} from a lateral (a) and posterior (b) viewpoints. b--c. RE between numerical approximations obtained with the AS and FS ($n=2$) methods. Results are presented from a lateral (c) and posterior (d) viewpoints.}
\label{fig_r2}
\end{figure*}

Figs.~\ref{fig_r2}a. and b. present the ratio $d/a$ for the model and sources described in Section~\ref{sec:real_1}. It can be clearly noted that there are some cortical regions for which the ratio is 0.5 or lower (dark blue in the figures), indicating that the FS ($n=2$) method may introduce unwanted errors compared to the AS approach. This assumption is confirmed in Figs.~\ref{fig_r2}c. and d., which show the $RE_s$ between both numerical approximations. It becomes evident that the use of the AS method outperforms the FS technique in these regions, improving up to 7\% in the resulting RE.  

\section{Discussion}\label{sec:8}

We presented analytical expressions for all the matrices and vectors involved in the FE formulation of the subtraction version of the EEG-FP considering dipolar sources. The adoption of these formulas allowed to minimise the errors for a given mesh discretisation, outperforming any subtraction-based approximation available in the literature. We found that the FS approach, considered the most accurate methodology for solving the EEG-FP based on the subtraction formulation~\cite{Beltrachini2019a}, introduces errors depending on the mesh refinement and source location through the ratio $d/a$. These inaccuracies were minimised with the AS technique, which prevents errors growing disproportionately no matter the mesh density and source position. Moreover, these advantages come at no extra computational cost, which resulted slightly lower than that corresponding to the FS method with $n=2$. This makes the AS technique the optimal option for solving the EEG-FP for any given model discretisation with regard to accuracy. 

The experiments presented provide new insights into the nature of the local errors due to the use of numerical quadrature formulas, as well as their impact in the overall solution of the EEG-FP. In the case of the surface element vector, $d$ (i.e. the distance between the source and the scalp surface) is always larger than 14mm. Then, from Fig.~\ref{fig_1}a., errors in $\bb^s$ introduced by a numerical scheme would become negligible for surface triangles with average side-length smaller than 14mm, which is always the case. This leaves the approximation of $\bb^v$ as the main responsible for the increased errors in the EEG-FP by means of the FS approach compared to those obtained with the AS method. In this case, errors will appear in the integration over elements with electrical conductivity different from that of the source compartment. We found that the use of numerical quadrature introduces errors depending on the ratio $d/a$ and the integration order~$n$ (Figs.~\ref{fig_1} and~\ref{fig_4}). These errors became noticeable for $d/a<0.5\ (0.25)$ in case of using $n=2\ (4)$, reaching values of $RE=1$ for $d/a\approx 0.05$. Such errors are avoided with the AS technique, which provides stable results regardless of $d/a$. 

Based on a real head model, we showed that the use of the AS method has clear benefits compared to the use of the FS approach, specially for sources located in the somatosensory and visual cortices, which are known for having an average thickness of approximately 1mm~\cite{Fischl2000}. Considering sources located in the midpoint between the GM/CSF and WM/GM surface layers, whose distance to the CSF and WM compartments is in the range 0.5-2.25mm~\cite{Fischl2000}, a tetrahedral mesh with side-lengths smaller than 1mm would be needed to avoid inaccuracies with the FS methodology (for $n=2$). Such refinement would need to take place in both CSF and WM regions, resulting in a highly-refined model and, consequently, in very high computational requirements. On the contrary, the AS method would not need such refinement to lead to even more accurate results, reducing both errors and computational demands. In the model used in Section~\ref{sec:real_1}, this resulted in the FS approach introducing an additional 7\% of RE in selected brain regions compared to the AS method, highlighting the value of the latter. 

Detailed head models are being increasingly considered for the analysis of both invasive and non-invasive EEG recordings. Recent results have demonstrated the value of incorporating the vessel network~\cite{Fiederer2016} and even the dura~\cite{Ramon2014} into the models used to localise sources of electrical activity. In this direction, evidence point to the incorporation of inhomogeneity of the electrical conductivity profile of the brain cortex, which presents significant laminar dependencies with radial/tangential anisotropies~\cite{Goto2010}. The addition of such level of detail would require a numerical method stable for extremely-eccentric sources located very close to conductivity jumps, being the AS approach the most appropriate (if not the only) option. 

Within all the methodologies available for simulating electromagnetic fields generated by brain current sources, the boundary element method is by far the most utilised. One of the reasons for this phenomenon is based on the great levels of accuracy reported for isotropic layered models, which are partly supported by the existence of analytical formulas for the potential integrals involved~\cite{deMunck1992b,Ferguson1994}. Such expressions were shown fundamental for speeding up simulations and avoiding errors due to the utilisation of a numerical quadrature scheme. We believe that the AS approach will have a similar effect with the FEM, extending its adoption and allowing optimal results in anisotropic and highly-detailed head representations.   

The MATLAB implementation of the method is publicly available through the FEMEG toolbox (\url{https://femeg.github.io}).

\section{Conclusions}~\label{sec:9}

We presented analytical expressions for all the integrals required by the FE formulation of the EEG-FP based on the subtraction approach. Such formulas were shown to increase the accuracy of the solution of the EEG-FP, preventing errors growing disproportionately independently of the mesh density and source position (unlike any other existing method). We demonstrated that the computational cost is similar to that obtained with the lowest order numerical scheme, making the AS method a competitive option in the field. We believe this approach will be considered as a gold standard technique for solving the EEG-FP considering fully-realistic head models.


\section*{Appendix 1: Proof of Proposition~\ref{prop_1}} \label{App_prop_1}

We use~\cite[A3.31]{VanBladel2007} to rewrite the surface divergence,
\[
\nabla_S \cdot \left( h(\bR) \frac{\bt}{t^2} \right) = \frac{\bt}{t^2}\cdot \nabla_Sh(\bR)+h(\br) \nabla_S\cdot\frac{\bt}{t^2}.
\]
The final result is obtained by noting that 
\begin{align*}
\nabla_S\cdot\frac{\bt}{t^2}=&\frac{\partial}{\partial u_a}\left( \frac{u_a}{u_a^2+v_a^2} \right)+ \frac{\partial}{\partial v_a}\left( \frac{u_a}{u_a^2+v_a^2} \right)\\
=&\frac{t^2-2u_a^2}{t^4}+\frac{t^2-2v_a^2}{t^4}=0.
\end{align*}

\section*{Appendix 2: Proof of Proposition~\ref{prop_2} } \label{App_prop_2}
For any two vector fields $\ba$ and $\bb$ in $\mathbb{R}^3$, it is easily verified that the following relation holds
\begin{equation}\label{eq:ap_prop_1}
\nabla(\ba \cdot \bb)=\nabla \ba^T \bb + \nabla \bb^T \ba,
\end{equation}
where $\nabla \ba$ is the gradient of the vector field with elements $\nabla\ba_{ij}=\partial a_i / \partial x_j$. Let $\ba=\bq$ and $\bb=\bR R^{-3}$. Direct computation of the partial derivatives of $\bR R^{-3}$ allows to check that its gradient is symmetric, i.e. $\nabla(\bR R^{-3})=\nabla(\bR R^{-3})^T$. Since $\bq$ is constant, the application of~\eqref{eq:ap_prop_1} yields
\begin{align*}
\hat{\bn} \cdot \nabla\left( \bq \cdot \bR R^{-3} \right)&=\hat{\bn}\cdot \left( \nabla \left(\bR R^{-3}\right)^T \bq  \right)= \bq^T \nabla \left(\bR R^{-3}\right) \hat{\bn}\\
&=\bq^T \nabla \left(\bR R^{-3}\right)^T \hat{\bn}=\bq \cdot \nabla \left( \hat{\bn} \cdot \bR R^{-3} \right).
\end{align*}

\section*{Appendix 3: Derivation of~\eqref{eq:I0_3} }\label{App1}

We operate as described in Section~\ref{sec:surf}, i.e. split the integral over $T$ as the sum of the integrals over $T_\epsilon$ and $T-T_\epsilon$, after which we take the limit $\epsilon \rightarrow 0$. Utilising~\eqref{eq:R5} we get
\begin{align}
\int_T \frac{dT}{R^5}&=\lim_{\epsilon \to 0} \left[ \int_{T_\epsilon} \frac{dT}{R^5} - \frac{1}{3}\int_{T-T_\epsilon} \!\!\!\!\!\!\!\nabla_S\cdot \left( \frac{1}{R^3} \frac{\bt}{t^2} \right) \right] \nonumber \\
& = \lim_{\epsilon \to 0} \left[ \int_{T_\epsilon} \frac{dT}{R^5} + \frac{1}{3}\int_{\partial T_\epsilon} \frac{\hat{\bm} \cdot \bt}{R^3 t^2} ds \right] - \frac{1}{3}\int_{\partial T} \frac{\hat{\bm} \cdot \bt}{R^3 t^2} ds. \label{eq:R5e}
\end{align}
Let $\alpha(\brho)$ be the portion of the circular arc $T_\epsilon$ lying within~$T$. This parameter is equal to $0$ if $\brho$ is outside $T$, $2\pi$ if $\brho$ is inside~$T$, $\pi$ if $\brho$ is on $\partial T$, or $0<\alpha<\pi$ if it is on a node. Utilising polar coordinates, the first two integrals yield
\begin{align*}
\int_{T_\epsilon} \frac{dT}{R^5}&=\int_0^\epsilon \int_0^{\alpha(\rho)}\! \! \! \! \! \!\frac{t}{\left( t^2 + w_0^2 \right)^{5/2}} d\phi dt\\
 &= -\frac{\alpha(\brho)}{3} \left( \frac{1}{(\epsilon^2 +w_0^2)^{3/2}} -\frac{1}{|w_0|^{3}} \right),
\end{align*}
and
\begin{align*}
\int_{\partial T_\epsilon} \frac{\hat{\bm} \cdot \bt}{R^3 t^2} ds= \int_0^{\alpha(\rho)}\! \! \! \! \! \!  \frac{d\phi}{(\epsilon^2 +w_0^2)^{3/2}}  =  \frac{\alpha(\brho)}{(\epsilon^2 +w_0^2)^{3/2}}. 
\end{align*}
The computation of the third integral in~\eqref{eq:R5e} is performed in~$\mathcal{L}_a$. Recalling the parameterisation $\bs(s)=\bt_i^0+s \hat{\bs}_i$, valid for $\bs\in\partial T$, and noting that $\hat{\bm}_i\cdot \bt=t_i^0$ for $\bt \in \partial T$, we get
\begin{align*}
\int_{\partial T} \frac{\hat{\bm} \cdot \bt}{R^3 t^2} ds=\sum_{i=1}^3 t_i^0 \int_{s_i^-}^{s_i^+} \! \! \! \frac{ds_i}{\left( (t_i^0)^2 + s_i^2 \right) \left( w_0^2 +(t_i^0)^2 + s_i^2 \right)^{3/2}} \\
= \sum_{i=1}^3  \frac{1}{|w_0|^3} \left( \tan^{-1}\frac{s_i^+ |w_0|}{t_i^0 R_i} - \tan^{-1}\frac{s_i^- |w_0|}{t_i^0 R_i}  \right)   -\frac{t_i^0}{w_0^2} R_i^s. 
\end{align*}
The final step consists of replacing these expressions into~\eqref{eq:R5e} and take the limit $\epsilon\rightarrow 0$. This leads to
\begin{align*}
\int_T \frac{dT}{R^5}=\frac{1}{3} \sum_{i=1}^3& \frac{t_i^0}{w_0^2} R_i^s + \frac{1}{3|w_0|^{3}} \Bigg[ \alpha(\brho)- \\ 
&\sum_{i=1}^3 \left( \tan^{-1}\frac{s_i^+ |w_0|}{t_i^0 R_i} - \tan^{-1}\frac{s_i^- |w_0|}{t_i^0 R_i}  \right) \Bigg].
\end{align*}
Noting that $\alpha(\brho)=\sum_{i=1}^3 \beta_i(w_0=0)$~\cite[eq.~(17)]{Graglia1993}, and employing the identity~\cite{Wilton1984}
\[
\tan^{-1} \frac{s_i^\pm}{t_i^0}-\tan^{-1}\frac{|w_0|s_i^\pm}{t_i^0 R_i^\pm} = \tan^{-1}\frac{t_i^0 s_i^\pm}{(R_i^0)^2+|w_0|R_i^\pm},
\]
we get the desired result.

\section*{Appendix 4: Derivation of~\eqref{eq:Iuava} }\label{sec:App3}

We multiply both sides of~\eqref{eq:I_dec} by $\{u_a,v_a\}^T$ and integrate over $T$,
\begin{align}
\int_T &
\begin{Bmatrix}
u_a \\ v_a
\end{Bmatrix}
\nabla \frac{w_0}{R^3} dT = w_0 \int_T 
\begin{Bmatrix}
u_a \\ v_a
\end{Bmatrix}
\nabla_S \frac{1}{R^3} dT \nonumber \\
&+ \hat{\bw} \int_T 
\begin{Bmatrix}
u_a \\ v_a
\end{Bmatrix}
\frac{dT}{R^3}-\hat{\bw} w_0^2 \int_T 
\begin{Bmatrix}
u_a \\ v_a
\end{Bmatrix}\frac{3 }{R^5} dT.\label{eq:app25}
\end{align}
After utilising~\eqref{eq:Rb}, the first term yields
\begin{align}
\int_T 
\begin{Bmatrix}
u_a \\ v_a
\end{Bmatrix}
 \nabla_S \frac{1}{R^3} dT =&  \int_T \nabla_S \left( 
\begin{Bmatrix}
u_a \\ v_a
\end{Bmatrix} 
 \frac{1}{R^3} \right) dT \nonumber\\ 
 &- 
\begin{Bmatrix}
 \hat{\bu} \\  \hat{\bv}
\end{Bmatrix} 
 \int_T\frac{1}{R^3} dT.\label{eq:ap21} 
\end{align}
The integrals in~\eqref{eq:ap21} are solved in $\mathcal{L}_a$. Using~\eqref{eq:gauss_2d} we get
\begin{align}
\int_T \nabla_S& \left( 
\begin{Bmatrix}
u_a \\ v_a
\end{Bmatrix} 
 \frac{1}{R^3} \right) dT = \sum_{i=1}^3 \hat{\bm}_i\int_{\partial T_i} 
 \begin{Bmatrix}
 u_a\\v_a
 \end{Bmatrix}
 \frac{ds_i}{R^3} \nonumber \\
&=\sum_{i=1}^3 \hat{\bm}_i   
\begin{Bmatrix}
 \hat{\bu} \\  \hat{\bv}
\end{Bmatrix}
\cdot \Bigg( \hat{\bs}_i \int_{s_i^-}^{s_i^+} \frac{s_i ds_i }{\left( w_0^2 + (t_i^0)^2 + s_i^2 \right)^{3/2}}\nonumber \\ 
&\hspace{3cm}+ \bt_i^0 \int_{s_i^-}^{s_i^+} \frac{ds_i}{\left( w_0^2 + (t_i^0)^2 + s_i^2 \right)^{3/2}} \Bigg) \nonumber\\
&=\sum_{i=1}^3 \hat{\bm}_i   
\begin{Bmatrix}
 \hat{\bu} \\  \hat{\bv}
\end{Bmatrix}
\cdot \left( \hat{\bs}_i R_i^d + \bt_i^0 R_i^s \right),\label{eq:ap22}
\end{align}
where we used 
\[
\begin{Bmatrix}
 u_a(s_i) \\  v_a(s_i)
\end{Bmatrix}
= (s_i \hat{\bs}_i + \bt_i^0) \cdot 
\begin{Bmatrix}
 \hat{\bu} \\  \hat{\bv}
\end{Bmatrix}.
\]
The second integral in~\eqref{eq:app25} is solved using~\eqref{eq:R3a}, resulting in
\begin{align}
&\int_T
\begin{Bmatrix}
u_a \\ v_a
\end{Bmatrix}
 \frac{1}{R^3} dT=-\int_T 
 \nabla_S\cdot \left(  
 \begin{Bmatrix}
\hat{\bu} \\ \hat{\bv}
\end{Bmatrix}
 \frac{1}{R} \right) dT \nonumber \\
 &\hspace{3mm}=-
 \begin{Bmatrix}
\hat{\bu} \\ \hat{\bv}
\end{Bmatrix}
\cdot \sum_{i=1}^3 \hat{\bm}_i \int_{\partial T_i} \frac{ds_i}{R}=
-
 \begin{Bmatrix}
\hat{\bu} \\ \hat{\bv}
\end{Bmatrix}
 \cdot \sum_{i=1}^3 \hat{\bm}_i f_{2i}.
\end{align}
Similarly, we apply~\eqref{eq:R5a} to solve the third integral in~\eqref{eq:app25}, obtaining
\begin{align}
&\int_T
\begin{Bmatrix}
u_a \\ v_a
\end{Bmatrix}
 \frac{3}{R^5} dT=-\int_T 
 \nabla_S\cdot \left(  
 \begin{Bmatrix}
\hat{\bu} \\ \hat{\bv}
\end{Bmatrix}
 \frac{1}{R^3} \right) dT \nonumber \\
 &\hspace{0mm}=-
 \begin{Bmatrix}
\hat{\bu} \\ \hat{\bv}
\end{Bmatrix}
\cdot \sum_{i=1}^3 \hat{\bm}_i \int_{\partial T_i} \frac{ds_i}{R^3}=
 -
 \begin{Bmatrix}
\hat{\bu} \\ \hat{\bv}
\end{Bmatrix} 
  \cdot \sum_{i=1}^3 \hat{\bm}_i R^s_i. \label{eq:ap24}
\end{align}
Finally, we replace~\eqref{eq:I0_2} and~\eqref{eq:ap22}--\eqref{eq:ap24} into~\eqref{eq:app25}, leading to the expected result.

\bibliographystyle{unsrt}
\bibliography{library}

\begin{thebibliography}{10}

\bibitem{Beltrachini2019b}
L.~Beltrachini.
\newblock {A finite element solution of the forward problem in EEG for
  multipolar sources}.
\newblock {\em IEEE Trans Neural Syst Rehab Eng.}, 2018.
\newblock in press.

\bibitem{deMunck1988a}
J.C. de~Munck, B.W. van Dijk, and H.~Spekreijse.
\newblock {Mathematical dipoles are adequate to describe realistic generators
  of human brain activity}.
\newblock {\em IEEE Trans Biomed Eng.}, 35(11):960--966, Nov. 1988.

\bibitem{Hamalainen1993}
M.~Hamalainen, R.~Hari, R.J. Ilmoniemi, J.~Knuutila, and O.V. Lounasmaa.
\newblock {Magnetoencephalography theory, instrumentation, and applications to
  noninvasive studies of the working human brain}.
\newblock {\em Rev Mod Phys.}, 65(2):413--497, 1993.

\bibitem{Baillet2001}
S.~Baillet, J.C. Mosher, and R.M. Leahy.
\newblock {Electromagnetic brain mapping}.
\newblock {\em IEEE Signal processing magazine}, 18(6):14--30, 2001.

\bibitem{deMunck1993}
J.C. de~Munck and M.J. Peters.
\newblock {A fast method to compute the potential in the multisphere model}.
\newblock {\em IEEE Trans Biomed Eng.}, 40(11):1166--1174, 1993.

\bibitem{vonEllenrieder2012}
N.~von Ellenrieder, L.~Beltrachini, and C.H. Muravchik.
\newblock {Electrode and brain modeling in stereo-EEG}.
\newblock {\em Clin. Neurophysiol.}, 123:1745--1754, 2012.

\bibitem{Cuffin1990}
BN~Cuffin.
\newblock {Effects of Head Shape on EEG’s and MEG’s}.
\newblock {\em IEEE Trans Biomed Eng.}, 37(1):44--52, 1990.

\bibitem{Huiskamp1999}
G.~Huiskamp, M.~Vroeijenstijn, R.~van Dijk, G.~Wieneke, and A.C. van Huffelen.
\newblock {The need for correct realistic geometry in the inverse EEG problem}.
\newblock {\em IEEE Trans Biomed Eng.}, 46(11):1281--1287, 1999.

\bibitem{Vatta2010}
F.~Vatta, F.~Meneghini, F.~Esposito, S.~Mininel, and F.~Di~Salle.
\newblock {Realistic and Spherical Head Modeling for EEG Forward Problem
  Solution: A Comparative Cortex-Based Analysis}.
\newblock {\em Comput Intell Neurosci.}, 2010:972060, 2010.

\bibitem{Beltrachini2011}
L.~Beltrachini, A.~Blenkmann, N.~von Ellenrieder, A.~Petroni, H.~Urquina,
  F.~Manes, A.~Ib\'{a}\~{n}ez, and C.H. Muravchik.
\newblock {Impact of head models in N170 component source imaging: results in
  control subjects and ADHD patients}.
\newblock {\em J Phys Conf Series}, 332(1):972060, 2011.

\bibitem{Wolters2006}
C.H. Wolters, A.~Anwander, X.~Tricoche, D.~Weinstein, M.A. Koch, and R.S.
  MacLeod.
\newblock {Influence of tissue conductivity anisotropy on {EEG/MEG} field and
  return current computation in a realistic head model: A simulation and
  visualization study using high-resolution finite element modeling}.
\newblock {\em NeuroImage}, 30:813--826, 2006.

\bibitem{vonEllenrieder2014b}
N.~von Ellenrieder, L.~Beltrachini, C.H. Muravchik, and J.~Gotman.
\newblock {Extent of cortical generators visible on the scalp: Effect of a
  subdural grid}.
\newblock {\em Neuroimage}, 101:787--795, 2014.

\bibitem{Beltrachini2019a}
L.~Beltrachini.
\newblock {Sensitivity of the projected subtraction approach to mesh
  degeneracies and its impact on the forward problem in EEG}.
\newblock {\em IEEE Trans Biomed Eng.}, 66(1):273--282, 2019.

\bibitem{Wolters2007}
C.H. Wolters, H.~K\"{o}stler, C.~M\"{o}ller, J.~H\"{a}rdtlein, L.~Grasedyck,
  and W.~Hackbusch.
\newblock {Numerical Mathematics of the Subtraction Method for the Modeling of
  a Current Dipole in EEG Source Reconstruction Using Finite Element Head
  Models}.
\newblock {\em SIAM J Sci Comput.}, 30(1):24--45, 2007.

\bibitem{Wolters2007c}
C.H. Wolters, H.~K\"ostler, C.~M\"oller, J.~H\"ardtlein, and A.~Anwander.
\newblock {Numerical approaches for dipole modeling in finite element method
  based source analysis}.
\newblock {\em Int Congress Series}, 1300:189--192, 2007.

\bibitem{Lew2009}
S.~Lew, C.H. Wolters, T.~Dierkes, C.~Röerc, and R.S. MacLeod.
\newblock {Accuracy and run-time comparison for different potential approaches
  and iterative solvers in finite element method based EEG source analysis}.
\newblock {\em Appl. Numer. Maths.}, 59:1970--1988, 2009.

\bibitem{Vorwerk2012}
J.~Vorwerk, M.~Clerc, M.~Burger, and C.H. Wolters.
\newblock {Comparison of Boundary Element and Finite Element Approaches to the
  EEG Forward Problem}.
\newblock {\em Biomed Tech Biomed Eng.}, 57:795--798, 2012.

\bibitem{Wilton1984}
D.R. Wilton, S.~Rao, A.~Glisson, D.~Shaubert, O.~Al-Bundak, and C.~Butler.
\newblock {Potential integrals for uniform and linear source distributions on
  polygonal and polyhedral domains}.
\newblock {\em IEEE Trans Antennas Propagation}, AP-32:276--1281, 1984.

\bibitem{Drechsler2009}
F.~Drechsler, C.H. Wolters, T.~Dierkes, H.~Si, and L.~Grasedyck.
\newblock {A full subtraction approach for finite element method based source
  analysis using constrained Delaunay tetrahedralisation}.
\newblock {\em Neuroimage}, 46:1055--1065, 2009.

\bibitem{Hutton2003}
D.V. Hutton.
\newblock {\em Fundamentals of Finite Element Analysis}.
\newblock McGraw-Hill Science, UK, 2003.

\bibitem{Beltrachini2015a}
L.~Beltrachini, Z.A. Taylor, and A.~Frangi.
\newblock {A parametric finite element solution of the generalised
  Bloch–Torrey equation for arbitrary domains}.
\newblock {\em J Magn Reson.}, 259:126--134, 2015.

\bibitem{Horn1991}
R.A. Horn and C.R. Johnson.
\newblock {\em Topics in Matrix Analysis}.
\newblock Cambridge University Press, UK, 1991.

\bibitem{Graglia1993}
R.D. Graglia.
\newblock {On the Numerical Integration of the Linear Shape Functions Times the
  3-D Green’s Function or its Gradient on a Plane Triangle}.
\newblock {\em IEEE Trans Antennas Propagation}, 41(10):1448--1455, 1993.

\bibitem{VanBladel2007}
J.~Van Bladel.
\newblock {\em {Electromagnetic Fields, second edition}}.
\newblock IEEE Press, John Wiley \& Sons, NJ, USA, 2007.

\bibitem{Zienkiewicz2013}
O.C. Zienkiewicz, R.L. Taylor, and J.Z. Zhu.
\newblock {\em The Finite Element Method: Its Basis and Fundamentals}.
\newblock Butterworth-Heinemann, UK, 2013.

\bibitem{Baumann1997}
S.B. Baumann, D.R. Wozny, S.K. Kelly, and F.M. Meno.
\newblock {The electrical conductivity of human cerebrospinal fluid at body
  temperature}.
\newblock {\em IEEE Transactions on Biomedical Engineering}, 44(3):20--23,
  1997.

\bibitem{Ramon2006}
C.~Ramon, P.H. Schimpf, and J.~Haueisen.
\newblock {Influence of head models on EEG simulations and inverse source
  localizations}.
\newblock {\em Biomed. Eng. Online}, 5:10, 2006.

\bibitem{Dannhauer2011}
M.~Dannhauer, B.~Lanfer, C.H. Wolters, and T.R. Kn\"{o}sche.
\newblock {Modeling of the human skull in {EEG} source analysis}.
\newblock {\em Human Brain Mapping}, 32:1383--1399, 2011.

\bibitem{McCann2019}
H.~McCann, G.~Pisano, and L.~Beltrachini.
\newblock {Variation in reported human head tissue electrical conductivity
  values}.
\newblock {\em bioRxiv}, 2019.
\newblock doi:10.1101/511006.

\bibitem{Fang2009}
Q.~Fang and D.~Boas.
\newblock {Tetrahedral mesh generation from volumetric binary and gray-scale
  images}.
\newblock In {\em Proc. IEEE Intl. Symp. Biomed. Imaging}, pages 1142--1145,
  2009.

\bibitem{Koay2011}
C.G. Koay.
\newblock {Analytically exact spiral scheme for generating uniformly
  distributed points on the unit sphere}.
\newblock {\em J Comput Sci.}, 2:88--91, 2011.

\bibitem{Aubert2006}
B.~Aubert-Broche, A.C. Evans, and D.L. Collins.
\newblock {A new improved version of the realistic digital brain phantom}.
\newblock {\em NeuroImage}, 32(1):138--145, 2006.

\bibitem{Dale1999}
A.~Dale, B.~Fischl, and M.I. Sereno.
\newblock {Cortical Surface-Based Analysis I: Segmentation and Surface
  Reconstruction}.
\newblock {\em Neuroimage}, 9:179--194, 1999.

\bibitem{Fischl2000}
B.~Fischl and A.~Dale.
\newblock {Measuring the thickness of the human cerebral cortex from magnetic
  resonance images}.
\newblock {\em Proc Natl Acad Sci USA.}, 97(20):11050--5, 2000.

\bibitem{Fiederer2016}
L.D.J. Fiederer, J.~Vorwerk, F.~Lucka, M.~Dannhauer, S.~Yang, M.~Dümpelmann,
  A.~Schulze-Bonhage~A. Aertsen, O.~Speck, C.H. Wolters, and T.~Ball.
\newblock {The role of blood vessels in high-resolution volume conductor head
  modeling of EEG}.
\newblock {\em Neuroimage}, 128:193--208, 2016.

\bibitem{Ramon2014}
C.~Ramon, P.~Garguilo, E.A. Fridgeirsson, and J.~Haueisen.
\newblock {Changes in scalp potentials and spatial smoothing effects of
  inclusion of dura layer in human head models for EEG simulations}.
\newblock {\em Front Neuroeng.}, 7:32, 2014.

\bibitem{Goto2010}
T.~Goto, R.~Hatanaka, T.~Ogawa, A.~Sumiyoshi, J.J. Riera, and R.~Kawashima.
\newblock {An evaluation of the conductivity profile in the somatosensory
  barrel cortex of Wistar rats}.
\newblock {\em J Neurophysiol.}, 104(6):3388--412, 2010.

\bibitem{deMunck1992b}
J.C. de~Munck.
\newblock {A linear discretization of the volume conductor boundary integral
  equation using analytically integrated elements}.
\newblock {\em IEEE Trans Biomed Eng.}, 39(9):986--90, 1992.

\bibitem{Ferguson1994}
A.S. Ferguson, X.~Zhang, and G.~Stroink.
\newblock {A complete linear discretisation for calculating the magnetic field
  using the boundary element method}.
\newblock {\em IEEE Trans Biomed Eng.}, 41(5):455--460, 1994.

\end{thebibliography}

\end{document}